\documentclass[12pt]{amsart} \hoffset -1cm \voffset -1in
\usepackage{marvosym}
\usepackage{amsmath}
\usepackage{amssymb}
\usepackage{amsthm}
\usepackage{color}
\usepackage{graphicx}
\usepackage[applemac]{inputenc}
\usepackage{phonetic}
\usepackage{pst-all}
\usepackage{pb-diagram}
\usepackage[T1]{fontenc}
\usepackage{stmaryrd}
\usepackage{url}

\newtheorem{prop}{Proposition}[section]
\newtheorem{lema}[prop]{Lemma}
\newtheorem{teor}[prop]{Theorem}
\newtheorem{hecho}[prop]{Fact}
\newtheorem{obs}[prop]{Observation}
\newtheorem{preg}[prop]{Question}
\newtheorem{corol}[prop]{Corollary}
\newtheorem{conj}[prop]{Conjecture}

\theoremstyle{definition}
\newtheorem{definicion}[prop]{Definition}
\newtheorem{ejem}[prop]{Example}
\newtheorem{rmk}[prop]{Remark}
\newtheorem{const}[prop]{Construction}

\def\Z{\mathbb Z}
\def\p{\varphi}
\def\tp{\mathrm{tp}}

\def\acl{\mathrm{acl}}
\def\dcl{\mathrm{dcl}}
\def\cb{\mathrm{Cb}}
\def\cupdot{\mathbin{\dot\cup}}
\def\bigcupdot{\mathop{\hbox to 0pt{\hskip.43em$\dot{\ }$\hss}\bigcup}}

\def\Ind#1#2{#1\setbox0=\hbox{$#1x$}\kern\wd0\hbox to 0pt{\hss$#1\mid$\hss}
\lower.9\ht0\hbox to 0pt{\hss$#1\smile$\hss}\kern\wd0}
\def\Notind#1#2{#1\setbox0=\hbox{$#1x$}\kern\wd0\hbox to 0pt{\mathchardef
\nn=12854\hss$#1\nn$\kern1.4\wd0\hss}\hbox to
0pt{\hss$#1\mid$\hss}\lower.9\ht0 \hbox to
0pt{\hss$#1\smile$\hss}\kern\wd0}

\def\ind{\mathop{\mathpalette\Ind{}}}                    
\def\nind{\mathop{\mathpalette\Notind{}}}          
\newcommand{\bp}{\begin{prop}}
\newcommand{\ep}{\end{prop}}
\newcommand{\bd}{\begin{definicion}}
\newcommand{\ed}{\end{definicion}}
\newcommand{\bej}{\begin{ejem}}
\newcommand{\eej}{\end{ejem}}
\newcommand{\bl}{\begin{lema}}
\newcommand{\el}{\end{lema}}
\newcommand{\bh}{\begin{hecho}}
\newcommand{\eh}{\end{hecho}}
\newcommand{\bpreg}{\begin{preg}}
\newcommand{\epreg}{\end{preg}}
\newcommand{\bo}{\begin{obs}}
\newcommand{\eo}{\end{obs}}
\newcommand{\bcon}{\begin{conj}}
\newcommand{\econ}{\end{conj}}
\newcommand{\brmk}{\begin{rmk}}
\newcommand{\ermk}{\end{rmk}}
\newcommand{\bc}{\begin{corol}}
\newcommand{\ec}{\end{corol}}
\newcommand{\bconst}{\begin{const}}
\newcommand{\econst}{\end{const}}
\newcommand{\bdem}{\begin{proof}}
\newcommand{\edem}{\end{proof}}
\newcommand{\benum}{\begin{enumerate}}
\newcommand{\eenum}{\end{enumerate}}
\newcommand{\bitem}{\begin{itemize}}
\newcommand{\eitem}{\end{itemize}}
\newcommand{\be}{\begin{ejem}}
\newcommand{\ee}{\end{ejem}}
\newcommand{\bt}{\begin{teor}}
\newcommand{\et}{\end{teor}}

\newcommand{\mult}{m}

\newcommand{\Le}{\mathcal{L}}

\title{Unimodularity Unified}
\author{Darío García and Frank O. Wagner}

\address{Univ Lyon, Universit\'e Claude Bernard Lyon 1, CNRS UMR 5208, Institut Camille Jordan, 43 blvd. du 11 novembre 1918, F-69622 Villeurbanne cedex, France}
\address{{\it Current address of the first author:} School of Mathematics, University of Leeds, LS2 9JT. Leeds, UK.}
\email{D.Garcia@leeds.ac.uk}
\email{wagner@math.univ-lyon1.fr}

\thanks{The first author is partially supported by the ``Programme Avenir Lyon Saint-Etiene'' (ANR-11-IDEX-0007) and the People Programme (Marie Curie Actions) of the European Union Seventh Framework Programme (FP7/2007-2013) under REA grant agreement PCOFUND-GA-2013-609102, through the PRESTIGE program coordinated by Campus France. The second author is partially supported by ValCoMo (ANR-13-BS01-0006).}
\date{03/10/2016}
\begin{document}

\begin{abstract}
Unimodularity is localized to a complete stationary type, and its properties are analysed. Some variants of unimodularity for definable and type-definable sets are introduced, and the relationship between these different notions is studied. In particular, it is shown that all notions coincide for non-multidimensional theories where the dimensions are associated to strongly minimal types.
\end{abstract}
\maketitle

\section*{Introduction}
Unimodularity was defined by Hrushovski in \cite{Hru} where he proved that a unimodular strongly minimal set is one-based, thus generalising Zilber's result that a locally finite strongly minimal set is 1-based. Recently, Hrushovski has re-visited unimodularity in the context of pseudofinite structures, aiming to develop an intersection theory for definable pseudofinite sets.

It was claimed in \cite{Hru} that unimodularity was equivalent to an {\em a priori} weaker notion called \emph{functional unimodularity} in \cite{Elw} and \cite{EJMR}. This was then used by Elwes as part of a proof that measurable stable structures are $1$-based \cite[Lemma 6.4]{Elw}, and was repeated in \cite{MS} and the survey article  \cite{EM}. In an attempt to clarify the situation, Pillay and Kestner \cite{Ke-Pi} have distinguished two types of functional unimodularity: one for definable sets and one for type-definable sets. They also studied the relationships between various notions and definitions, mainly in the context of strongly minimal structures. In particular, they showed that for strongly minimal theories, unimodularity is equivalent to functional unimodularity for arbitrary types, and is also equivalent to the structures being measurable in the sense of \cite{MS2}. They also presented an example intended to be a strongly minimal set which is functionally unimodular but not unimodular. However, the example actually turns out not to be functionally unimodular; in fact our Theorem \ref{nonmulti} states that all variants of unimodularity coincide for non-multidimensional theories where the dimensions are associated to strongly minimal types. 

This paper can be seen as yet another attempt to clarify the situation, and is organized as follows: In Section 1 we introduce the notion of a uniform correspondence, measurability of a (partial) type, and commensurability between (partial) types, and develop the basic properties. In Section 2 we introduce the concept of correspondence unimodularity and functional unimodularity for complete types, partial types and definable sets, and give a correction to Proposition 3.2 in \cite{Ke-Pi}. The main result in this section is Theorem \ref{unimodequiv}, which states that unimodularity is equivalent to both correspondence unimodularity and to functional unimodularity for complete types, and Theorem \ref{ct-t}, which says that in an $\omega$-stable theory unimodularity is equivalent to both correspondence and functional unimodularity for partial types.

In Section 3 we localize unimodularity to complete stationary types, and finally show that all concepts coincide for non-multidimensional theories where the dimensions are associated to strongly minimal types, and in particular for $\aleph_1$-categorical theories and groups of finite Morley rank. 

We use standard model-theoretic notation and work in some big sufficiently saturated and ultrahomogeneous monster model of the theory. Lower case letters $a,b,c,$ etc. will denote finite tuples. If a tuple $a$ is algebraic over $b$, we use $m(a/b)$ for the (finite) number of realizations of $\tp(a/b)$. We shall not distinguish between singletons and tuples, or between real and imaginary elements (i.e. we work in $T^{eq}$).

%

\section{Correspondences}
\bd[Correspondence] Let $\pi$ and $\pi'$ be two type-definable sets. \benum
\item A {\em correspondence} between $\pi$ and $\pi'$ is a non-empty type-definable set $C(x,y)\vdash\pi(x)\times\pi'(y)$ such that all fibres $C_x=\{y\models\pi':C(x,y)\}$ and $C^y=\{x\models\pi:C(x,y)\}$ are finite. If $\pi'=\pi$ we call $C$ a correspondence {\em on} $\pi$.
\item A correspondence $C$ is {\em complete} if it is a complete type. 
\item A correspondence $C$ is {\em uniform} if the fibre sizes $k_C=|C_x|$ and $\ell_C=|C^y|$ are constant, independently of $x\models\pi$ and $y\models\pi'$.
\item A \emph{$(k,\ell)$-correspondence} is a uniform correspondence with $k=k_C$ and $\ell=\ell_C$.
\item For a uniform correspondence $C$, the {\em ratio} of $C$ is $m_C=\frac{k_C}{\ell_C}$.
\item A correspondence $C$ is {\em balanced} if it is uniform and $k_C=\ell_C$ (equivalently,  $m_C=1$).
\eenum
If $\pi$, $\pi'$ and $C$ are all type-definable over some parameters $A$, we say that $C$ is {\em over} $A$.
\ed

Note that a uniform correspondence is actually relatively definable, by compactness. If $C(x,y)$ is a correspondence between $\pi(x)$ and $\pi'(y)$, then $C^{-1}(y,x)=C(x,y)$ is a correspondence between $\pi'(y)$ and $\pi(x)$. Clearly, $(C^{-1})_y=C^y$ and $(C^{-1})^x=C_x$. So $C^{-1}$ is uniform/complete/balanced if and only if $C$ is.\medskip

Correspondences between complete types are particularly well-behaved.
\bl \label{fmcorresp} Let $C$ be a correspondence between a complete type $p$ and some partial type $\pi(y)$, all over the same parameters $A$. Then:
\benum
\item $|C_x|$ does not depend on $x\models p$.
\item $C$ can be written as the disjoint union of finitely many complete correspondences $C=C_0\cupdot \cdots \cupdot C_n$, with $n\leq |C_x|$.
\eenum
\el
\bdem \benum
\item If $a,a'\models p$, then there is an automorphism $\sigma$ fixing $A$ with $\sigma(a)=a'$. Then $C_{a'}=\sigma(C_a)$, so $|C_{a'}|=|C_a|$.
\item If $\tp(a_i,b_i/A)$ for $i\in I$ are the completions of $C$, then $a_i\models p$ and we may assume $a_i=a_0$ for all $i\in I$. But then $b_i\in C_{a_0}$; since the types $\tp(a_0,b_i/A)$ are all different, we have $b_i\not=b_j$ for $i\not=j$, and $|I|\le|C_{a_0}|$. It follows that $C=\bigcupdot_{i\in I}C_i$ with $C_i=\tp(a_i,b_i/A)$.\qedhere
\eenum
\edem

\bc \label{corfmcorresp} A correspondence $C$ between complete types is automatically uniform, and if all its completions have the same ratio $m$, then $m_C=m$.
\ec
\bdem Suppose $C(x,y)$ is a correspondence between complete types $p(x)$ and $q(y)$. Then $|C_x|=k_C$ and $|C^y|=|(C^{-1})_y|=\ell_C$ are constant for $x\models p$ and $y\models q$ by Lemma \ref{fmcorresp}, hence the correspondence is uniform. If $C_0,\ldots,C_n$ are the completions of $C(x,y)$, then
\begin{equation}\label{sumcomp}
k=|C_x|=\Big|\bigcupdot_{i=0}^n (C_i)_x\Big|=\sum_{i=0}^n k_{C_i} \text{\ \ and\ \ } \ell=|C^y|=\Big|\bigcupdot_{i=0}^n (C_i)^y\Big|=\sum_{i=0}^n \ell_{C_i}.
\end{equation}
If all the completions $C_i$ have the same ratio $m$, then $k_{C_i}=m\ell_{C_i}$ for all $i$, whence $k=m\ell$ and $m_C=m$.
\edem 

\bd[Measurable, Commensurable]\label{mesdefn} Let $\pi$ be a partial type over $A$. We say that $\pi$ is {\em measurable over $A$} if every $A$-type-definable uniform correspondence $C$ on $\pi$ is balanced.\par
Two partial types $\pi$ and $\pi'$ over $A$ are {\em commensurable over $A$} if there is a uniform correspondence $C$ from $\pi$ to $\pi'$, and for any other uniform correspondence $C'$ over $A$ between $\pi$ and $\pi'$ one has $m_{C'}=m_C$. In this case we put $m_\pi^{\pi'}=m_C$.
If $\pi$ is measurable over any $B\supseteq A$, we say that $\pi$ is measurable; if $\pi$ and $\pi'$ are commensurable over any $B\supseteq A$ we say that they are commensurable.
\ed
Thus $\pi$ is measurable (over $A$) if and only if $\pi$ and $\pi$ are commensurable (over $A$). It follows from Corollary \ref{corfmcorresp} that for complete types we may restrict ourselves to complete correspondences in Definition \ref{mesdefn}.

If $B\supseteq A$ and $\pi$ and $\pi'$ are commensurable over $B$, and if there is a correspondence between $\pi$ and $\pi'$ over $A$, then $\pi$ and $\pi'$ are commensurable over $A$. However, commensurability or measurability over $A$ need not imply commensurability or measurability over $B$.

\bl\label{commem-compl} Two complete types $p$ and $q$ are commensurable over $A$ if and only if there is a complete correspondence $C$ over $A$ between $p$ and $q$, and all such complete correspondences take the same value $m_C=m_p^q$.\el

\bdem The left to right direction follows directly from the definitions. Conversely, let $C_0,\ldots,C_n$ be the completions of $C$. By (\ref{sumcomp}),
$$k_C=\sum_{i=0}^n k_{C_i}=\sum_{i=0}^n m_p^q\cdot \ell_{C_i}=m_p^q \ell_C.$$
This yields the result.\edem

We shall now study composition of correspondences. 
\bd[Composition] Let $\pi$, $\pi'$ and $\pi''$ be partial types over $A$, and suppose $C$, $C'$ are correspondences between $\pi$ and $\pi'$ and between $\pi'$ and $\pi''$, respectively. The composition $C'\circ C$ is defined by
\[(a,c)\in C'\circ C \Leftrightarrow \exists b\,[(a,b)\in C\land(b,c)\in C'].\]\ed
By compactness and saturation, $C'\circ C$ is type-definable; note that any witness $b$ for the existential quantifier must automatically satisfy $\pi'$. It is clear that $(C'\circ C)_a$ and $(C'\circ C)^c$ are finite for every $a\models\pi$ and $c\models\pi''$, so $C'\circ C$ is a correspondence between $\pi$ and $\pi''$. 

If $\pi$ and $\pi''$ are complete types over $A$, then $C'\circ C$ can be written as a finite union $D_0\cup \cdots \cup D_n$ of complete correspondences between $\pi$ and $\pi''$ by Lemma \ref{fmcorresp}, each of which is uniform by Corollary \ref{corfmcorresp}. If moreover $C$ and $C'$ are both uniform (for instance if $\pi'$ is also complete), given $(a,c)\in D_i$ define 
$$r_i=|\{b\models\pi': (a,b)\in C \text{  and  }(b,c)\in C' \}|.$$
Since $D_i$ is complete, this number only depends on $D_i$ and not on the choice of $(a,c)\models D_i$. Then for $a\models\pi$
\begin{equation} \label{sumcompos}\begin{aligned}
k_C\cdot k_{C'}&=|\{(b,c):(a,b)\in C\land (b,c)\in C'\}|\\
&=|\bigcup_{i\le n}\{(b,c):(a,c)\in D_i\land (a,b)\in C\land (b,c)\in C'\}|=
\sum_{i=0}^n r_i\cdot k_{D_i}.\end{aligned}\end{equation}
Similarly $\displaystyle{\ell_C\cdot \ell_{C'}=\sum_{i=0}^n r_i\cdot \ell_{D_i}}$.

\bp \label{composition} Let $p$, $q$ and $r$ be complete types, and suppose $C$ is a correspondence between $p$ and $q$ and $C'$ is a correspondence between $q$ and $r$, all over $A$. If $p$ and $r$ are commensurable over $A$, then $m_{C'\circ C}=m_C\cdot m_{C'}$.
\ep
\bdem By Lemma \ref{fmcorresp} the correspondences $C$, $C'$ and $C'\circ C$ are all uniform; let $(D_i:i\le n)$ be the finitely many completions of $C'\circ C$. Since $p$ and $r$ are commensurable over $A$, we have that $m_{D_i}=m_p^r$ for every $i\leq n$. By (\ref{sumcompos}) we obtain
\[k_{C}\cdot k_{C'}=\sum_{i=1}^n r_i\cdot k_{D_i}=\sum_{i=1}^n r_i\cdot (m_p^r\cdot \ell_{D_i}) = m_p^r\sum_{i=1}^n r_i\cdot \ell_{D_i} = m_p^r\cdot (\ell_C\cdot \ell_{C'}),\]
whence
\[m_{C'\circ C}=m_p^r=\dfrac{k_{C}\cdot k_{C'}}{\ell_{C}\cdot \ell_{C'}}=m_C\cdot m_{C'}\]
\edem

\bc\label{measure-commem} Let $p$ and $q$ be complete types over $A$.\benum
\item Suppose there is a correspondence $C$ between $p$ and $q$. If $p$ is measurable over $A$, then so is $q$, and $p$ and $q$ are commensurable over $A$.
\item If $p$ and $q$ are commensurable over $A$, then $p$ and $q$ are both measurable over $A$.
\item For any three complete commensurable types $p$, $q$ and $r$ over $A$ we have
$m_p^q\,m_q^r=m_p^r$.\eenum\ec

\bdem\benum\item If $C'$ is any other correspondence between $p$ and $q$ over $A$, then $C'^{-1}(y,x)=C'(x,y)$ is a correspondence from $q$ to $p$. Clearly $m_{C'^{-1}}=m_{C'}^{-1}$. By Lemma
\ref{composition} we have 
$$1=m_p^p=m_C\cdot m_{C'^{-1}}=m_C/m_{C'},$$
so $m_{C'}=m_C=m_p^q$. Hence $p$ and $q$ are commensurable over $A$.
\item Suppose that $p$ and $q$ are commensurable over $A$. If $C$ is a complete correspondence on $p$ over $A$, then $m_Cm_p^q=m_p^q$ by Proposition \ref{composition}, and $m_C=1$. Thus $p$ is measurable over $A$; measurability of $q$ over $A$ follows by symmetry.
\item This follows immediately from Proposition \ref{composition}.\eenum\edem

\bt\label{ct-meas=>t-meas}Let $\pi$ be a partial type over $A$ and suppose $MR(\pi)<\infty$. If all completions of $\pi$ over $A$ of maximal Morley rank are measurable over $A$, so is $\pi$.\et
\bdem Suppose $C$ is a $(k_C,\ell_C)$-correspondence on $\pi$ over $A$. Let $(p_i:i\in I)$ be the finitely many completions of $\pi$ over $A$ of maximal Morley rank. Then for all $i,j\in I$, if $C_{ij}=C\cap(p_i\times p_j)$ is non-empty, it is a correspondence between $p_i$ and $p_j$, so the two types are commensurable by Corollary \ref{measure-commem}. If $C_{ij}=\emptyset$ put $k_{C_{ij}}=\ell_{C_{ij}}=0$. Put
$$I_0=\{i\in I:p_1\text{ and $p_i$ are commensurable over $A$}\}.$$
If $(a,b)\in C$ with $a\models p_i$ for some $i\in I_0$, then by interalgebraicity
$$RM(b/A)=RM(ab/A)=RM(a/A),$$
so $b\models p_j$ for some $j\in I_0$. It follows that for each $i\in I_0$
$$\sum_{j\in I_0}k_{C_{ij}}=k_C\qquad\text{and}\qquad{\sum_{j\in I_0} \ell_{C_{ji}}=\ell_C}.$$
For $i\in I_0$ put $m_i=m_{p_1}^{p_i}$. If $C_{ij}\not=\emptyset$ we have $m_j=m_i\cdot m_{C_{ij}}$ by Corollary \ref{measure-commem}, that is $$m_j\cdot \ell_{C_{ij}}=m_i\cdot k_{C_{ij}}.$$
Note that the latter equation trivially holds if $C_{ij}=\emptyset$.

Put $\mu=\sum_{i\in I_0} m_i$. Then $\mu\not=0$ and
\begin{align*}
\mu\cdot k_C &=\sum_{i\in I_0} (m_i\cdot k_C)=\sum_{i\in I_0}\big(m_i \sum_{j\in I_0} k_{C_{ij}}\big)=\sum_{i\in I_0}\sum_{j\in I_0} (m_i\cdot k_{C_{ij}})\\
&=\sum_{i\in I_0}\sum_{j\in I_0} (m_j\cdot \ell_{C_{ij}})=\sum_{j\in I_0}\sum_{i\in I_0} (m_j\cdot \ell_{C_{ij}})=\sum_{j\in I_0} \big(m_j\sum_{i\in I_0}\ell_{C_{ij}}\big)\\
&=\sum_{j\in I_0} (m_j\cdot \ell_C)=\mu\cdot \ell_C.
\end{align*}
It follows that $k_C=\ell_C$.\edem

\bej  Let $M=\Z\times 2^\omega$ in the language $\{f,E_n:n\in\omega\}$, where the $E_n$ are equivalence relations with $2^n$ classes given by
$$(z,\eta)E_n(z',\eta')\Leftrightarrow z\equiv z'\mod2^n$$
and 
$$f(z,\eta)=(z+1,\eta\circ S),$$
where $S$ is the successor function on $\omega$. Then $E_{n+1}$ cuts each $E_n$-class in half, and $f:M\to M$ is a surjective function with fibres of size two. Moreover, $x E_n y\Leftrightarrow f(x) E_n f(y)$, and for any $m\in M$ the $2^n$ elements $m, f(m),f^2(m),\ldots,f^{2^n-1}(m)$ are in different $E_n$-classes. This theory is complete of Lascar rank one, but not $\omega$-stable. Every stationary complete type is measurable, but the model itself (equivalently, the partial type $x=x$) is not. So $\omega$-stability is necessary in Theorem \ref{ct-meas=>t-meas}.\eej

\section{Unimodularity and its variations.}
We shall now study the relationship between unimodularity introduced in \cite{Hru}, functional unimodularity and its variants formally introduced in \cite{Ke-Pi}, and \emph{correspondence unimodularity} for definable sets, complete types or types. We start with some definitions.

\bd [Unimodularity] A complete theory is \emph{unimodular} if for any two tuples $a$, $b$ and parameters $A$ in the monster model, if $a\equiv_Ab$ and $a$ and $b$ are interalgebraic over $A$, then $m(a/Ab)=m(b/Aa)$.
\ed
\bl\label{uniequiv} A theory is unimodular if and only if every complete type is measurable over its domain.\el
\bdem Let $p(x)$ be a complete type. Note that two realizations $a,b\models p$ are $A$-inter\-al\-ge\-braic if and only if $C=\tp(a,b/A)$ is a complete correspondence on $p$ over $A$. Then $m(b/Aa)=k_C$ and $m(a/Ab)=\ell_C$. So $m(a/Ab)=m(b/Aa)$ if and only if $C$ is balanced. By Corollary \ref{corfmcorresp}, any correspondence on $p$ is balanced if and only if all complete correspondences on $p$ are balanced. Thus shows the equivalence.\edem

\bd [Functional unimodularity] Let $T$ be a complete theory. Then $T$ is
\benum
\item \emph{functionally unimodular} (FU) if for any two definable sets $X$ and $Y$ we have:\begin{itemize}
\item[(*)] If two definable functions $f,g:X\to Y$ have constant fibre sizes $k$ and $\ell$ respectively, then $k=\ell$;\end{itemize}
\item \emph{functionally unimodular for types} (FU-t) if property (*) holds for any type-defin\-able sets $X$, $Y$;
\item \emph{functionally unimodular for complete types} (FU-ct) if property (*) holds for any complete types $X$, $Y$.
\eenum
\ed

Kestner and Pillay  \cite{Ke-Pi} proved that if $T$ is strongly minimal, then unimodularity is equivalent to functional unimodularity for types, and in this case it is also equivalent to \emph{MS-measurability}. We shall now show that functional unimodularity allows finitely many exceptional finite fibres.
\bp\label{finiteexceptions} Let $X$ and $Y$ be two infinite
definable sets, and $f,g:X\to Y$ two definable functions with finite 
fibres, such that $|f^{-1}(y)|=k$ and $|g^{-1}(y)|=\ell$ for all but finitely 
many $y\in Y$. If $k\not=\ell$, there are definable sets $X'$ and $Y'$, as well as definable functions $f',g':X'\to Y'$ such that the fibres of $f'$ and $g'$ have constant sizes $k$ and $\ell$, respectively.\ep
\bdem Put
$$Y_0=\{y\in Y:|f^{-1}(y)|\not=k\text{ or }|g^{-1}(y)|\not=\ell\}.$$
Let $F=f^{-1}(Y_0)$ and $G=g^{-1}(Y_0)$. Without loss of generality we may 
assume that $|F|\le|G|$; modifying $f$ definably on finitely many points we 
may further assume $F\subseteq G$. Put 
$$X''=X\setminus F,\quad Y''=Y\setminus Y_0,\quad G'=G\setminus 
F,\quad\text{and}\quad f''=f\restriction_{X''}:X''\to Y''.$$
Then $f''$ has constant fibre size $k$, and 
$$g\restriction_{X''\setminus G'}:X''\setminus G'\to Y'$$
has constant fibre size $\ell$. Put $n=|G'|$.

{\em Case 1: $k<\ell$.} Let $n'=\ell-k$. Let $P$ be a set of cardinality $kn$ 
and $Q$ a set of cardinality $n$. Put 
$$X'=(X\times n')\cup P,\quad Y'=(Y''\times n')\cup Q,$$
and define $f':X'\to Y'$ via $f'((y,i))=(f''(y),i)$ for $(y,i)\in X''\times n'$, and $f':P\to Q$ 
arbitrarily with fibres of constant size $k$. Finally, define $g':X'\to Y'$ 
via $g'((y,i))=(g(y),i)$ for $(y,i)\in(X''\setminus G')\times n'$, and 
$g':(G'\times n')\cup P\to Q$ arbitrarily with fibres of constant size $\ell$, 
which is possible since 
$$|(G'\times n')\cup P|=nn'+kn=n(\ell-k+k)=\ell n=\ell\,|Q|.$$

{\em Case 2: $\ell<k$.} Let $n'=k-\ell-1$. Let $Q\subset Y''$ have cardinality 
$n$, and put $P=f''^{-1}(Q)\subset X''$, of cardinality $kn$. We choose $Q$ such 
that $P\cap G'=\emptyset$. Put 
$$X'=(X''\times n')\cup((X''\setminus P)\times\{n'\}),\quad 
Y'=(Y''\times n')\cup((Y''\setminus Q)\times\{n'\}),$$
and define $f':X'\to Y'$ via $f'((y,i))=(f''(y),i)$, with fibres of 
constant size $k$. Note that the map
$$g'':(X''\setminus G')\times(n'+1)\to Y''\times(n'+1)$$
defined by $g''((y,i))=(g(y),i)$ has constant fibre size $\ell$. Now $X'$ has 
$$|P|-|G'\times(n'+1)|=kn-n(n'+1)=(k-(k-\ell))n=\ell n$$ points less than
$(X''\setminus G')\times(n'+1)$, and $Y'$ has $|Q|=n$ points less than 
$Y''\times(n'+1)$. Modifying $g''$ on finitely many points, we can thus define 
a map $g':X'\to Y'$ with constant fibre size $\ell$.\edem
\bc\label{finexc} Let $T$ be functionally unimodular. If $X$ and $Y$ are two definable sets, and $f,g:X\to Y$ are two definable maps of constant fibre sizes $k$ and $\ell$, respectively, except for finitely many exceptional fibres which are still finite, then $k=\ell$.\ec
\bdem This follows immediately from Proposition \ref{finiteexceptions}.\edem

\bej\label{CUct-CUt} Consider the structure $M=\langle 2^{<\omega},S\rangle$ where is $S$ is interpreted as the successor relation, that is, $D\models S(a,b)$ if and only if $a\text{\textasciicircum}0=b$ or $a\text{\textasciicircum}1=b$. This structure is strongly minimal, and was proposed in \cite{Ke-Pi} as an example of a strongly minimal structure which is functionally unimodular but not unimodular. The non-unimodularity follows from the fact that if $S(a,b)$ holds, then $a$ and $b$ are interalgebraic but $\mult(a/b)=1\neq 2=\mult(b/a)$.

Contrary to \cite[Proposition 3.2]{Ke-Pi}, in fact this structure is {\em not} functionally unimodular: The identity function id$_M$ is clearly $1$-to-$1$, while the predecessor function $f$ defined by the formula 
\[\varphi(x,y)=S(y,x) \lor (\forall z(\neg S(z,x))\land x=y)\] is $2$-to-$1$ almost everywhere, with an exceptional fibre of size $3$ at $\emptyset$. So $M$ is not functionally unimodular by Corollary \ref{finexc}. This can also be seen directly: Add an additional point $\infty$ to the structure, and define $f'(x)=f(x)$ for $x\not=\emptyset$, and $f'(\emptyset)=f'(\infty)=\infty$. Then $f'$ is surjective and $2$-to-$1$ on $M\cup\{\infty\}$, contradicting functional unimodularity.\eej

\bd [Correspondence unimodularity] A complete theory $T$ is \emph{correspondence unimodular} (CU) if for any two definable sets $X$ and $Y$ we have:
\begin{itemize}
\item[(**)] If $C_1$ and $C_2$ are uniform correspondences between $X$ and $Y$, then $m_{C_1}=m_{C_2}$.
\end{itemize}
We say that $T$ is correspondence unimodular for (complete) types (CU-t and CU-ct, respectively), if (**) holds whenever $X$ and $Y$ are (complete) types.\ed

\bl \label{corr-itself} A theory $T$ is correspondence unimodular (resp.\ for types or complete types) if and only if all definable sets (resp.\ types or complete types) are measurable.
\el
\bdem $(\Rightarrow)$ Suppose $C$ is a uniform correspondence on $\pi$. Then $C^{-1}$ is again a uniform correspondence on $\pi$. By correspondence unimodularity, $m_C=m_{C^{-1}}=1/m_C$, whence $m_C=1$ and $C$ is balanced.

$(\Leftarrow)$ Suppose $C_1$, $C_2$ are uniform correspondences between $\pi_1$ and $\pi_2$. Define $C$ on $\pi_1\times\pi_2$ by
 \[(a_1,b_1)C(a_2,b_2) \Leftrightarrow a_1C_1b_2 \wedge a_2C_2b_1.\]
It is easy to see that $C$ is a uniform correspondence on $\pi_1\times\pi_2$, with
$$k_C=k_{C_1}\cdot \ell_{C_2}\qquad\text{and}\qquad\ell_C=k_{C_2}\cdot \ell_{C_1}.$$
By assumption $k_C=\ell_C$, whence $m_{C_1}=m_{C_2}$. So $T$ is correspondence unimodular.
\edem 

\bej \label{psfcorresp} It is easy to show that all pseudofinite structures are correspondence unimodular (for definable sets): If $M=\prod_{\mathcal{U}} M_i$ is an ultraproduct of finite structures and $C$ is a uniform correspondence on a definable set $X\subseteq M$, then in the finite structures $M_i$ we have that 
\[|C_i|=\big|\bigcup_{x\in X_i} \{(a,b)\in C_i:a=x\}\big|=\sum_{x\in X_i} |(C_i)_x|=|X_i|\cdot k_C\]
for $\mathcal{U}$-almost all indices $i$. Similarly, $|C_i|=|X_i|\cdot\ell_C$, whence $k_C=\ell_C$ and $C$ is balanced. Therefore all definable sets are measurable; by Lemma \ref{corr-itself} we have correspondence unimodularity.
\eej

We shall now identify various implications between the different notions of unimodularity. It is clear that functional unimodularity for types implies both functional unimodularity for complete types and for definable sets, and similarly for correspondence unimodularity. We shall show the implications given by the dotted arrows in the diagram below, sometimes under additionnal model-theoretic hypotheses.

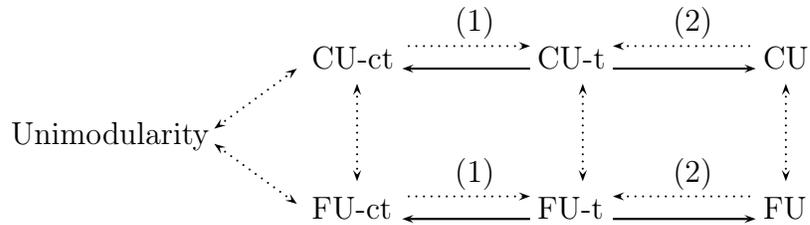
\begin{center}
\begin{pspicture}(2,0)(3,3)
\put(-3,1){Unimodularity}
\put(1,2){CU-ct}
\put(1,0){FU-ct}
\put(4,2){CU-t}
\put(4,0){FU-t}
\put(7,2){CU}
\put(7,0){FU}
\put(2.9,2.5){(1)}
\put(2.9,0.5){(1)}
\put(5.8,2.5){(2)}
\put(5.8,0.5){(2)}

\psline[linestyle=dotted,dotsep=2pt]{<->}(-0.3,1.2)(0.8,2)
\psline[linestyle=dotted,dotsep=2pt]{<->}(-0.3,1)(0.8,0.2)

\psline[linestyle=dotted,dotsep=2pt]{<->}(1.6,1.8)(1.6,0.5)
\psline[linestyle=dotted,dotsep=2pt]{<->}(4.6,1.8)(4.6,0.5)
\psline[linestyle=dotted,dotsep=2pt]{<->}(7.3,1.8)(7.3,0.5)

\psline{->}(3.9,0)(2.2,0)
\psline{->}(3.9,2)(2.2,2)

\psline[linestyle=dotted,dotsep=2pt]{->}(2.2,2.3)(3.9,2.3)
\psline[linestyle=dotted,dotsep=2pt]{->}(2.2,0.3)(3.9,0.3)

\psline{->}(5,2)(6.9,2)
\psline{->}(5,0)(6.9,0)

\psline[linestyle=dotted,dotsep=2pt]{->}(6.9,2.3)(5,2.3)
\psline[linestyle=dotted,dotsep=2pt]{->}(6.9,0.3)(5,0.3)
\end{pspicture}
\end{center}

\benum
\item $T$ $\omega$-stable
\item $T$ non-multidimensional, with strongly minimal dimensions\eenum

We first note that the functional and correspondence versions of unimodularity are equivalent.
\bp\label{C-F} A theory is functionally unimodular (resp.\ FU-t or FU-ct) if and only it is correspondence unimodular (resp.\ CU-t or CU-ct).\ep
\bdem $(\Rightarrow):$ Let $C(x,y)$ be a uniform correspondence on a definable set $X$ (resp.\ type-definable set or complete type). Note that if $X$ is a complete type, by Corollary \ref{corfmcorresp} we may assume that $C$ is complete. Consider the two functions $f,g:C\to X$, where $f$ is the projection to the first and $g$ the projection to the second coordinate. Then $f$ is $k_C$-to-$1$ and $g$ is $\ell_C$-to-$1$. By functional unimodularity (resp.\ FU-t or FU-ct) we have $k_C=\ell_C$, and $C$ is balanced. By Lemma \ref{corr-itself} we are done.

$(\Leftarrow):$ Suppose $X$ and $Y$ are type-definable sets, and $f,g:X\to Y$ are relatively definable surjective functions that are respectively $k$-to-$1$ and $\ell$-to-$1$. Consider the correspondence $C$ on $X$ defined by 
\[(a,a')\in C \Leftrightarrow f(a)=g(a').\]
Then $C$ is a $(\ell,k)$-correspondence on $X$, and $k=\ell$ by correspondence unimodularity.\edem
As a corollary, we obtain in general the equivalence between unimodularity and functional unimodularity for complete types, originally shown by Kestner and Pillay for strongly minimal theories.
\bt\label{unimodequiv} Let $T$ be a complete theory. The following are equivalent:\benum
\item $T$ is unimodular.
\item $T$ is correspondence unimodular for complete types.
\item $T$ is functionally unimodular for complete types.\eenum\et
\bdem This follows from Lemmas \ref{uniequiv} and \ref{corr-itself} and Proposition \ref{C-F}.\edem

Example \ref{CUct-CUt} shows that our next theorem does need $\omega$-stability.
\bt\label{ct-t} Let $T$ be $\omega$-stable unimodular. Then $T$ is correspondence unimodular for types.\et
\bdem This follows from Lemmas \ref{uniequiv} and \ref{corr-itself} and Theorem \ref{ct-meas=>t-meas}.\edem

The following is an example of a functionally unimodular structure which is not unimodular. We shall show in Theorem \ref{nonmulti} that for a non-multidimensional theory with strongly minimal dimensions, functional unimodularity does imply unimodularity.
\bej For each $n<\omega$, let $M_n=2^{<n}$. We consider $M_n$ as a finite structure in the language $\Le=\{R_i:i<\omega\}\cup \{f\}$ by interpreting the predicates as $R_i^{M_n}=\{\eta\in M_n: \operatorname{length}(\eta)=n-i\}$ for $i\leq n$, and $R_i^{M_n}=\emptyset$ for $i>n$. To interpret the function $f$ we put:
\[f(\eta))=\begin{cases}\eta\upharpoonright_{\operatorname{length(\eta)-1}} &\text{if $\operatorname{length}(\eta)>1$}\\ \emptyset &\text{if $\eta=\emptyset$}\end{cases}\]
Let $M=\prod_{\mathcal{U}}M_n$, where $\mathcal{U}$ is a non-principal ultrafilter over $\omega$. Note that in the ultraproduct, $f:M\to M$ is a definable function such that $f\upharpoonright_{R_i}:R_i\twoheadrightarrow R_{i+1}$ is a $2$-to-$1$ function.

Since $M$ is pseudofinite, it is correspondence unimodular (Example \ref{psfcorresp}). It is easy to check that $M$ is $\omega$-stable, even non-multidimensional of Morley rank $2$. However, $M$ is not correspondence unimodular for complete types: Consider the complete type given by $$q(x)=\{\neg R_i(x):i<\omega\}\cup\{f^i(x)\not=x:i<\omega\}.$$
Then $f(q)=q$, and $f\upharpoonright_q$ is $2$-to$1$, so $q$ is not measurable.
\eej

\section{Unimodularity for types}
Throughout this section we shall work  in a stable theory with elimination of imaginaries.
We first introduce some notions from geometric stability theory. For further reading, the reader can consult \cite{P} or \cite{Wa00}.
\bd Let $\pi$ be a partial type over $A$, and $\Sigma$ an $A$-invariant family of partial types. Then $\pi$ is\begin{itemize}
\item ({\em almost}) {\em $\Sigma$-internal} if for every realization $a$
of $\pi$ there is $B\ind_Aa$ and a tuple $\bar b$ of realizations of types in $\Sigma$
based on $B$, such that $a\in\dcl(B\bar b)$ (or $a\in\acl(B\bar b)$,
respectively).
\item {\em $\Sigma$-analysable} if for any realization $a$ of $\pi$ there are
$(a_i:i<\alpha)\in\dcl(Aa)$ such that $\tp(a_i/A,a_j:j<i)$ is
$\Sigma$-internal for all $i<\alpha$, and
$a\in\acl(A,a_i:i<\alpha)$. We call $\alpha$ the {\em length} of the analysis.\end{itemize}
\ed
We shall say that $a$ is (almost) $\Sigma$-internal or $\Sigma$-analysable over $b$ if $\tp(a/b)$ is.
\bd Two types $p\in S(A)$ and $q\in S(B)$ are {\em orthogonal} if for all $C\supseteq AB$, $a\models p$, and $b\models q$ with $a\ind_A C$ and $b\ind_B C$ we have $a\ind_C b$.\par
A type $p$ is {\em regular} if it is orthogonal to all its forking extensions.\par
A theory is {\em non-multidimensional} if every type is non-orthogonal to a type over $\emptyset$.\ed
Equivalently, a theory is non-multidimensional if there are only boundedly many pairwise orthogonal types.
\bd[Unimodularity] A complete stationary  type $p$ is {\em unimodular} if over any set $A$ of
parameters containing dom$(p)$, whenever $a$ and $b$ are $A$-interalgebraic realizations of the non-forking extension of $p$ to $A$, then $m(a/Ab)=m(b/Aa)$.\ed
\brmk Equivalently, $p$ is unimodular if all its non-forking extensions are measurable over their domain.\ermk
\bl\label{product} Let $p$ and $p'$ be unimodular stationary types of finite Lascar
rank over $A$. Let $aa'$ and $bb'$ be $A$-interalgebraic realizations of the free product
$p\otimes p'$. Suppose $a\ind_A b'$ and $a'\ind_A b$. Then $m(aa'/Abb')=m(bb'/Aaa')$.\el
\bdem By stationarity and independence, $a$ and $b$ both realize $p|Aa'$. Moreover, $b\in\acl(Aaa')$. By the Lascar equalities in finite rank, 
$$U(a/Aa'b)=U(aa'b/A)-U(a'b/A)=U(aa'b/A)-U(aa'/A)=U(b/Aaa')=0.$$
So $a$ and $b$ are $Aa'$-interalgebraic, whence $m(a/Aa'b)=m(b/Aaa')$ by unimodularity of $p$. Thus
$$m(bb'/Aaa')=m(b'/Aaa'b)\cdot m(b/Aaa')=m(b'/Aaa'b)\cdot m(a/Aba')=m(ab'/Aba').$$
Similarly, $b'$ and $a'$ are $Ab$-interalgebraic realizations of $p'|Ab$. So
$m(b'/Aba')=m(a'/Abb')$ by unimodularity of $p'$, and
$$m(ab'/Aba')=m(a/Aba'b')\cdot m(b'/Aba')=m(a/Aba'b')\cdot m(a'/Abb')=m(aa'/Abb').\qedhere$$\edem
\bc\label{orthogonal} If $p$ and $q$ are orthogonal unimodular stationary types of
finite Lascar rank, then their free product $p\otimes q$ is unimodular.\ec
\bdem This follows immediately from the definitions and Lemma \ref{product}.\edem
\bc\label{power} If $p$ is a unimodular regular stationary type of finite Lascar
rank, then the free power $p^{(n)}$ is unimodular for all $n\ge 1$.\ec
\bdem We can assume $p\in S(\emptyset)$. If $(a_i:i<n)$ and $(b_i:i<n)$ are two interalgebraic realizations of
$p^{(n)}$, put $\bar a=(a_i:i>0)$. Let $\tilde b=(b_i:b_i\nind\bar a)$. Since
$a_0\ind\bar a$ we have $a_0\ind\tilde b$. Let $\bar b\supseteq\tilde b$ be maximal
with $a_0\ind\bar b$. Then $\bar b$ has length $n-1$, and there is a unique
$b_j\notin\bar b$. Note that $b_j\ind\bar a$. As $\bar a$ and $\bar b$ satisfy
$p^{(n-1)}=:p'$, and $a_0$ and $b_j$ satisfy $p$, the hypotheses of Lemma
\ref{product} are satisfied, and we conclude.\edem

\bl\label{mlem} Let $\pi$ and $\pi'$ be partial types over $A$, and $A\subseteq B$. Put
$$\bar\pi(x):=\pi(x)\land x\ind_AB\qquad\text{and}\qquad\bar\pi'(y):=\pi'(y)\land y\ind_AB.$$If $C$ is a uniform correspondence between $\pi$ and $\pi'$ over $A$, then $C'=C\cap(\bar\pi\times \bar\pi')$ is a uniform correspondence between $\bar\pi$ and $\bar\pi'$ with $m_{C'}=m_C$.\el
\bdem For $a\models\bar\pi$ we have $a\ind_AB$. If $(a,b)\in C$, then $b\models\pi'$ and $b\in\acl(Aa)$, whence $b\ind_AB$ and $b\models\bar\pi'$. Thus $(a,b)\in C'$, and $|(C')_a|=k_C$. Similarly $|(C')^b|=\ell_C$ for all $b\models\bar\pi'$. Therefore $C'$ is uniform, with $k_{C'}=k_C$ and $\ell_{C'}=\ell_C$, whence $m_{C'}=m_C$.\edem
\bc\label{nf-unimod} Suppose $q$ is a non-forking extension of a stationary type $p$. Then $p$ is unimodular if and only if $q$ is unimodular.\ec
\bdem $(\Rightarrow)$ follows from the definition. For the converse, consider a non-forking extension $p'$ of $p$, and the common non-forking extension $q'$ of $p'\cup q$. Take $\pi=\pi'=p'$ and $\bar\pi=\bar\pi'=q'$ in Lemma \ref{mlem}. As $m_{C'}=1$ by measurability of $q$, we get $m_C=1$ and $p'$ is measurable. Hence $p$ is unimodular.\edem
\bc\label{interalg} Let $p$ and $q$ be stationary types over $A$ whose realizations are $A$-inter\-al\-ge\-braic. Suppose $p$ is unimodular.\benum
\item Then $q$ is unimodular, and $p$ and $q$ are commensurable.
\item If $p'$ and $q'$ are non-forking extension of $p$ and $q$ to the same domain, then $p'$ and $q'$ are again commensurable, and $m_p^q=m_{p'}^{q'}$.\eenum\ec
\bdem As $p$ is measurable, $p$ and $q$ are commensurable by Corollary \ref{measure-commem}. Moreover, $p'$ and $q'$ are also commensurable by Lemma \ref{mlem}, and $m_p^q=m_{p'}^{q'}$. Hence all non-forking extensions of $q$ are measurable, and $q$ is unimodular.\edem

\bc\label{internal} Let $P$ be an $\emptyset$-invariant family of unimodular weakly
minimal stationary types. If $q$ is almost $P$-internal, then $q$ is
unimodular.\ec
\bdem Since $q$ is almost $P$-internal, there is a realization $a\models q$, some set $A$ of parameters independent of $a$, and realisations $\bar b$ of types in $P$ over $A$ with $a\in\acl(A\bar b)$. As $P$ consists of weakly minimal types, we may assume that $\bar b$ is independent over $A$. Let $\bar b=\bar b'\bar b''$, where $\bar b'$ is a maximal subtuple of $\bar b$ independent of $a$ over $A$. Then $\tp(a/A\bar b')$ is a non-forking extension of $q$, and $a$ and $\bar b''$ are interalgebraic over $A\bar b'$ by weak minimality of the types in $P$. Moreover, $\bar b''$ is independent over $A\bar b'$. The result now
follows from Corollaries  \ref{orthogonal}, \ref{power}, \ref{nf-unimod} and~\ref{interalg}.\edem

We now turn to analysability. Let us first consider an example which shows that non-multidimensionality is necessary in Theorem \ref{nonmulti}.
\bej Let $E$ be an equivalence relation with infinitely many infinite classes, and
$f$ a unary surjective function with fibres of size two, such that $xEx'\Leftrightarrow
f(x)Ef(x')$ and that neither $f$ nor the induced relation $f_E$ on $E$-classes have
any non-oriented cycles (and in particular $\neg xEf(x)$). It is easy to see that this theory is multidimensional of
Morley rank $2$; one dimension is carried by the type $\tp(a_E)$ of the $E$-classes,
and the other dimensions by $\tp(a/a_E)$, for any $a$. Each dimension has Morley
rank $1$ and is unimodular. Nevertheless, $\tp(a)$ is clearly not unimodular, as
$a\equiv f(a)$, $m(f(a)/a)=1$ but $m(a/f(a))=2$.\eej

\bt\label{analysable} Let $P$ be a set of unimodular strongly minimal
types over $\emptyset$. Then any $P$-analysable stationary type is unimodular.\et
\bdem By Corollary \ref{nf-unimod} we may add parameters to the language and suppose that the types in $P$ are over $\emptyset$. Note that as the types in $P$ are strongly minimal, any $P$-analysable stationary type $q$ is contained in a definable set $\p$ which is $P$-analysable of finite length. Then $\p$ is non-multidimensional, and its dimensions are strongly minimal. So $\p$ is $\omega$-stable by \cite[Corollaire 2.14]{Po}.

We shall use induction on the length of a $P$-analysis of $q$. If it is $1$, then $q$ is almost $P$-internal, and we are done by Corollary \ref{internal}. 

So suppose $q$ has a $P$-analysis of length $n+1$. For $b\models q$ put
$$B=\{e\in\acl(b):\tp(e)\text{ has a $P$-analysis of length at most $n$}\},$$
the $n$-th $P$-level $\ell_n^P(b)$ (see \cite[Definition 3.1]{PW}).
Put $A=B\cap\dcl(b)$. If $e\in B$ and $e'\equiv_be$, then $e'\in B$, and there are only finitely many such $e'$. Let $\bar e$ be any imaginary element coding this finite set. Then $\bar e\in\dcl(b)$, and $\bar e\in\dcl\{e':e'\equiv_be\}\subseteq B$, so $\bar e\in A$. Hence $B=\acl(A)$. Moreover, the type $\tp(b/A)$ is stationary, as $\tp(b/B)$ is stationary, $b\ind_AB$, and for every $A$-definable finite equivalence relation $E$ the class $b_E$ of $b$ modulo $E$ is in $$\dcl(Ab)\cap\acl(B)=\dcl(b)\cap B=A.$$
By $\omega$-stability of $\p$ we can choose $a\in A$ such that $b\ind_aA$ and $\tp(b/a)$ is stationary; note that then $A=\dcl(a)$. Since $\tp(b)$ has a $P$-analysis of length $n+1$, the type $\tp(b/B)$ and thus also $\tp(b/a)$ is almost $P$-internal, whence unimodular. Finally, $\tp(a)$ is stationary since $\tp(b)$ is, and unimodular by inductive hypothesis.

If $b'\models q$ and $b$ and $b'$ are interalgebraic, choose $a'$ with $a'b'\equiv ab$.
Note that $\cb(a'/b)$ is definable over a Morley sequence in $\tp(a/b')$, and thus has a $P$-analysis of length at most $n$. It follows that $\cb(a'/b)\in B$ and $a'\ind_B b$, whence $a'\ind_a b$. Similarly, $a\ind_{a'}b'$. But 
$$a\in\dcl(b)\subseteq\acl(b)=\acl(b')\qquad\mbox{and}\qquad a'\in\dcl(b')\subseteq\acl(b')=\acl(b),$$
so the independences above imply that $a$ and $a'$ are interalgebraic.

By stationarity of $\tp(b/a)$, the independence $b\ind_aa'$ and unimodularity of $\tp(a)$ we have 
$$m(a'/ab)=m(a'/a)=m(a/a')=m(a/a'b').$$

Since $\tp(b)$ is almost $P$-internal, there is $D\ind_ab$ containing $a$ and some tuple $d$ of realizations of types in $P$ over $D$ such that $b$ and $d$ are $D$-interalgebraic. As $\tp(d)$ is $P$-internal, it is unimodular by Corollary \ref{internal}, as is $\tp(b/a)$. Put $p=\tp(b/a)$, $q=\tp(b'/a')$ and $r=\tp(d)$. Let $p^*$ and $q^*$ be the non-forking extensions of $p$ and $q$ to $aa'$. As $b$ and $b'$ are $aa'$-interalgebraic and $p$ is unimodular, $p^*$ and $q^*$ are commensurable and
$$m_{p^*}^{q^*}=\dfrac{m(b'/aa'b)}{m(b/aa'b')}.$$
Let $\sigma$ be a strong $\emptyset$-automorphism mapping $a$ to $a'$, and put
$D'=\sigma(D)$. Let $p'$, and $r'$ be the non-forking extensions of $p$ and $r$ to $D$, and $q'$ and $r^*$ the non-forking extensions of $q$ and $r$ to $D'$. As $p$ is unimodular, $p'$ and $q'$ are commensurable, as are $\sigma(p')=q'$ and $\sigma(r')=r^*$. Clearly $m_{p'}^{r'}=m_{q'}^{r^*}$.

Finally, let $p''$, $q''$ and $r''$ be the non-forking extensions of $p$, $q$ and $r$ to $DD'$. Then $p''$, $q''$ and $r''$ are commensurable by Corollary \ref{interalg}, and by Lemma \ref{composition} we get
$$m_{p''}^{r''}=m_{p''}^{q''}m_{q''}^{r''}.$$
But now by Corollary \ref{interalg} again,
$$\frac{m(b'/aa'b)}{m(b/aa'b')}=m_{p^*}^{q^*}=m_{p''}^{q''}=\frac{m_{p''}^{r''}}{m_{q''}^{r''}}=\frac{m_{p'}^{r'}}{m_{q'}^{r^*}}=1.$$
Hence $m(b'/aa'b)=m(b/aa'b')$. As $a\in\dcl(b)$ and $a'\in\dcl(b')$ we finally obtain
$$\begin{aligned}m(b/b')&=m(ab/a'b')=m(b/aa'b')m(a/a'b')\\
&=m(b'/aa'b)m(a'/ab)=m(a'b'/ab)=m(b'/b).\end{aligned}$$
It follows that $q$ is unimodular.\edem
\bt\label{nonmulti} Let $T$ be a non-multidimensional theory whose dimensions are associated to strongly minimal types. The following are equivalent:\benum
\item $T$ is unimodular.
\item $T$ is functionally unimodular.
\item All strongly minimal types are unimodular.\eenum\et
\bdem $(1)\Rightarrow(2):$ By \cite[Corollaire 2.14]{Po} the theory $T$ is $\omega$-stable, so unimodularity implies functional unimodularity for partial types by Theorem \ref{ct-t}. Functional unimodularity (for sets) follows.

$(2)\Rightarrow(3):$ Let $p$ be a strongly minimal type which is not unimodular. We may assume $p$ is over $\emptyset$. So there are interalgebraic realizations $a,b\models p$ with $m(a/b)\not= m(b/a)$. Then $\tp(a,b)$ has Morley rank $1$. Choose definable sets $X\in\tp(a,b)$ and $Y\in p$ of Morley rank $1$, such that $Y$ has Morley degree $1$ and $X\subset Y\times Y$. Consider the functions $f,g:X\to Y$, where $f$ is the projection to the first coordinate, and $g$ is the projection to the second coordinate. Restricting $Y$ we may assume that $f$ has fibres of size at most $m(b/Aa)$, and $g$ has fibres of size at most $m(a/Ab)$. As $Y$ is strongly minimal and the fibre sizes are bounded, there are only a finite number of exceptional fibres, of size less than $m(b/a)$ for $f$ and of size less than $m(a/b)$ for $g$. By Proposition \ref{finiteexceptions} there are definable sets $X'$ and $Y'$ and definable functions $f',g':X'\to Y'$ whose fibres all have size $m(b/a)$ and $m(a/b)$, respectively. As $m(b/a)\not=m(a/b)$, this contradicts functional unimodularity.

$(3)\Rightarrow(1):$ Let $P$ be a set of strongly minimal types containing a representative for each dimension. Then every type is $P$-analysable, and hence unimodular by Theorem \ref{analysable}.\edem
Examples of non-multidimensional theories whose dimensions are associated to strongly minimal types are almost strongly minimal theories, uncountably categorical theories, and groups of finite Morley rank.

\section{Further remarks}
Although we have defined unimodularity for arbitrary stationary types, we could only show that it is well-behaved for types of finite rank. The problem obviously comes from the fact that in infinite rank, say close to a regular type $p$, we should work with $p$-closure rather than algebraic closure, which is unbounded. Thus multiplicity is not the correct measure.

A possibility might be to define {\em Lascar unimodularity}: Let us say that a stationary type $p$ over $A$ is {\em Lascar unimodular} if for any realizations $a,b\models p$ we have $U(a/Ab)=U(b/Aa)$. Theories of finite Lascar rank are clearly Lascar unimodular. This notion may be particularly pertinent if $p$ is a regular type, as then $a$ and $b$ are dependent if and only if either one is in the $p$-closure of the other. However, we have not studied the properties of Lascar unimodularity, nor have we looked for interesting examples.

Another question concerns unimodularity for non-stationary types. Section $1$ of our paper does not assume stationarity, so one might be tempted to develop unimodularity, at least for Lascar strong types, in a simple theory in analogy with Section $3$.

\end{document}